\documentclass[12pt,a4paper]{article}
\usepackage[utf8]{inputenc}
\usepackage{amssymb,amsmath,theorem}
\usepackage{hyperref}

\newcommand\qed{\ifhmode\unskip\nobreak\fi\quad            %  квадрат в конце
   \ifmmode\square\else\hbox{$\square$}\fi}                %   доказательств
\newcommand\pin{\kern.0833em}                              %  половина тонкого пробела
\let\eps=\varepsilon \let\kappa=\varkappa                  %  греческие буквы

\theorempreskipamount=\medskipamount
\theorempostskipamount=\medskipamount
\newcommand\proof[1]{\textit{Proof#1\pin}}

\newtheorem{theorem}{Theorem}%[section]
\newtheorem{lemma}[theorem]{Lemma}
\newtheorem{corollary}[theorem]{Corollary}

\DeclareMathOperator{\var}{var}
\DeclareMathOperator{\diam}{diam}

\begin{document}

\begin{flushright}
\itshape Dedicated to V.\,G. Krotov
\end{flushright}

\begin{center}
{\Large\bfseries Criterion for the absolute continuity of curves\\[4pt] in metric spaces}

\bigskip

\itshape V.\,I. Bakhtin

\medskip

Belarusian State University, Minsk, Belarus (e-mail: bakhtin@tut.by) %ORCID 0000-0001-5782-5378
\end{center}

\vspace{-24pt}

\renewcommand{\abstractname}{}
\begin{abstract} \noindent
It is proved that a parameterized curve in a metric space $X$ is absolutely continuous if and only if its composition with any Lipschitz function on $X$ is absolutely continuous.
\end{abstract}

\smallskip

%\quad\parbox{13.9cm}
\textbf{Keywords:} {\itshape absolutely continuous curve, metric space, Lipschitz function}

\smallskip

\textbf{2020 MSC:} 51F30, 53C23

\bigskip\medskip

This article was written solely due to the brilliant report by A.\,I. Tyulenev presented on April 29, 2024 at a seminar led by V.\,G. Krotov at the Belarusian State University. There was proved a criterion for the absolute continuity of curves in metric spaces satisfying the doubling condition (which automatically implies finiteness of the Hausdorff dimension of those spaces).

While thinking about the necessity of the above mentioned doubling condition I was very surprised to discover that it could simply be omitted (see Theorem \ref{..1} below). Had been informed about this, A.\,I. Tyulenev replied that he had also discovered this. Moreover, in collaboration with R.\,D. Oleynik, he proved a criterion for $p\pin$-absolute continuity of a curve for any $p\in [1,+\infty]$ (in this article this is done only for $p=1$), and a criterion for a curve to belong to the generalized Sobolev class~$W^{1,p}$. In essence, Theorem \ref{..1} below is just a special case of Theorem 1.1 from \cite{Tyulenev}. However, its proof is much simpler compared to \cite{Tyulenev} and is carried out using elementary reasoning.

%Already after his report, A.\,I. Tyulenev, in collaboration with R.\,D. Oleynik, independently and almost simultaneously, were also able to get rid of the doubling condition. Moreover, in \cite{Tyulenev} they proved a criterion for $p\pin$-absolute continuity of a curve for any $p\in [1,+\infty]$ (in this article this is done only for $p=1$), and a criterion for a curve to belong to the generalized Sobolev class~$W^{1,p}$. In essence, Theorem \ref{..1} below is just a special case of Theorem 1.1 from \cite{Tyulenev}. However, its proof is much simpler compared to \cite{Tyulenev} and is carried out using elementary reasoning.

\smallskip

A mapping $\gamma\colon [a,b] \to X$ of a segment $[a,b]\subset\mathbb R$ to a metric space $(X,\rho)$ is called \emph{abso\-lut\-ely continuous} (notation $\gamma \in \mathit{AC}([a,b], X)$, %see \cite[p. 108]{Royden}, \cite[p. 252]{Nielsen}, \cite[p. 128]{Athreya} for the case of
see \cite[p.~128]{Athreya}, \cite[p.~252]{Nielsen} for the case of a real-valued function) if for every $\eps>0$ there exists $\delta>0$ such that for any finite collection of pairwise disjoint intervals $(a_i,b_i)\subset [a,b]$ the following implication is valid:
\begin{equation*}
 \sum_i (b_i-a_i) < \delta \quad \Longrightarrow \quad \sum_i \rho(\gamma(a_i),\gamma(b_i)) <\eps.
\end{equation*}

\begin{theorem} \label{..1}
A continuous mapping\/ $\gamma\colon [a,b] \to X$ is absolutely continuous if and only if for any Lipschitz function\/ $h\colon X\to \mathbb R$ the composition\/ $h\circ \gamma$ is absolutely continuous.
\end{theorem}

%The necessary part of this theorem is obvious. To prove its sufficient part, we need the following three lemmas.

The necessary part of this theorem is obvious. To prove its sufficient part we need three lemmas. The following Lemma \ref{..2} is a generalization of the well-known McShane theorem \cite[Theorem 1]{McShane} on the extension of a Lipschitz function (which follows immediately from Lemma~\ref{..2} for $L' =L$).

\begin{lemma} \label{..2}
Suppose a real-valued function\/ $h$ defined on a certain subset\/ $X'\subset X$ satisfies Lipschitz condition with constant\/ $L'$$:$
\begin{equation*}
 \big|h(x') -h(x'')\big| \le L' \rho(x',x''), \qquad x'\!\pin,\pin x''\in X'.
\end{equation*}
Then for each\/ $L\ge L'$ the functions
\begin{gather} \label{,,1}
  h^+_L(x) =\inf_{x'\in X'}\{\pin h(x') +L\rho(x',x)\pin\}, \qquad x\in X, \\[3pt]
  h^-_L(x) =\sup_{x'\in X'}\{\pin h(x') -L\rho(x',x)\pin\}, \qquad x\in X \label{,,2}
\end{gather}
have the following properties\/$:$

\smallskip

a) their restrictions to\/ $X'$ coincide with\/ $h$$;$

\smallskip

b) they are\/ $L$-Lipschitz $($satisfy Lipschitz condition with constant\/ $L$$);$

\smallskip

c) for all\/ $x\in X$ we have
\begin{equation}\label{,,3}
  h^+_L(x) - h^-_L(x) \ge 2(L-L')\pin\rho(x,X').
\end{equation}
\end{lemma}

\proof. If $x\in X'$ then the infimum in \eqref{,,1} is attained at $x' =x$ and coincides with $h(x)$. Similarly, the supremum in \eqref{,,2} is attained at $x' =x$ and coincides with $h(x)$. This proves~a).

Property b) is obvious. To prove c) consider the difference
\begin{equation} \label{,,4}
  h^+_L(x) -h^-_L(x) \pin=\pin
  \inf_{x'\!,\, x''\in X'}\bigl\{ h(x') -h(x'') +L\rho(x',x) +L\rho(x'',x)\bigr\}.
\end{equation}
Note that in the case $\rho(x',x'')\le 2\rho(x,X')$ we have
\begin{equation*}
   h(x')-h(x'') +L\rho(x',x) +L\rho(x'',x) \ge -L'\rho(x',x'') +2L\rho(x,X') \ge 2(L-L')\pin\rho(x,X'),
\end{equation*}
and in the case $\rho(x',x'')\ge 2\rho(x,X')$ we have
\begin{equation*}
   h(x')-h(x'') +L\rho(x',x) +L\rho(x'',x) \ge -L'\rho(x',x'') +L\rho(x',x'') \ge 2(L-L')\pin\rho(x,X').
\end{equation*}
In both cases, these inequalities and \eqref{,,4} imply \eqref{,,3}. \qed

\begin{corollary} \label{..3}
Under the conditions of Lemma \ref{..2} the function\/ $h$ extends to an $L$-Lipschitz function on\/ $X'\cup\{x\}$ when given any value\/ $h(x)$ between\/ $h^-_L(x)$ and\/ $h^+_L(x)$.
\end{corollary}

\begin{lemma} \label{..4}
Suppose\/ $h$ is an\/ $L'$-Lipschitz function defined on a certain subset\/ $X'\subset X$, and\/ $\gamma\colon [a,b]\to X$ is a continuous mapping satisfying the condition\/ $\rho(\gamma[a,b],X') >0$. Then for any constant\/ $L>L'$ there exist a partition\/ $T=\{t_0,t_1,\dots,t_n\}$ of the segment\/ $[a,b]$ $($i.\,e., $a=t_0<t_1<\dotsc <t_n=b$$)$ and an\/ $L$-Lipschitz extension of\/ $h$ to\/ $X'\cup \gamma(T)$ such that
\begin{equation}\label{,,5}
  \var(h\circ\gamma,T) \pin:=\pin
  \sum_{i=1}^n \big| h(\gamma(t_{i-1})) -h(\gamma(t_i))\big| \pin\ge\pin L\rho(\gamma(a),\gamma(b)).
\end{equation}
\end{lemma}

\proof. Let $h^{\pm}_L(x)$ be the functions defined in Lemma \ref{..2}. Then, by point c) of this lemma, $h^{+}_L(x) > h^{-}_L(x)$ for all $x\in \gamma([a,b])$.

Starting from $t_0 =a$, we will sequentially choose $t_i$ and define alternating values
\begin{equation}\label{,,6}
  h(\gamma(t_{2i})) =h^-_L(\gamma(t_{2i})), \qquad   h(\gamma(t_{2i+1})) =h^+_L(\gamma(t_{2i+1})),
\end{equation}
except the last point $t_n=b$, where it will be
\begin{equation}\label{,,7}
  h^-_L(\gamma(t_{n})) \pin\le\pin h(\gamma(t_{n})) \pin\le\pin  h^+_L(\gamma(t_{n})),
\end{equation}
according to the next rule: for even $i$ we take
\begin{gather} \label{,,8}
  t_{i+1} \pin=\pin \sup\bigl\{ t\in [t_i,b] : h^-_L(\gamma(t_i)) +L\rho(\gamma(t_i),\gamma(t)) \le h^+_L(\gamma(t))\bigr\}, \\[6pt] \label{,,9}
  h(\gamma(t_{i+1})) \pin=\pin h^-_L(\gamma(t_i)) +L\rho(\gamma(t_i),\gamma(t_{i+1})),
\end{gather}
and for odd $i$ we take
\begin{gather}  \label{,,10}
  t_{i+1} \pin=\pin \sup\bigl\{ t\in [t_i,b] : h^+_L(\gamma(t_i)) - L\rho(\gamma(t_i),\gamma(t)) \ge h^-_L(\gamma(t))\bigr\}, \\[6pt]  \label{,,11}
  h(\gamma(t_{i+1})) \pin=\pin h^+_L(\gamma(t_i)) -L\rho(\gamma(t_i),\gamma(t_{i+1})).
\end{gather}
Since the mapping $\gamma\colon [a,b]\to X$ is uniformly continuous, and the functions $h^{+}_L(x)$ and $h^{-}_L(x)$ are $L$-Lipschitz (see point b) of Lemma \ref{..2}), this procedure will reach $t_n=b$ in a finite number of steps.

Now \eqref{,,6}, \eqref{,,9}, \eqref{,,11} imply the equalities
\begin{equation} \label{,,12}
  \bigl|h(\gamma(t_{i-1})) -h(\gamma(t_{i}))\bigr| \pin=\pin L\rho(\gamma(t_{i-1}),\gamma(t_i)), \qquad i=1,\dots,n.
\end{equation}
Summing them up, we obtain \eqref{,,5}.

\pagebreak[1]

Let us check Lipschitz condition for the restriction of $h$ to $\gamma(T)$. For any pair of adjacent points $\gamma(t_{i-1})$, $\gamma(t_{i})$ it is proved in \eqref{,,12}. Consider an arbitrary pair of points $\gamma(t_{i})$, $\gamma(t_{j})$, where $j\ge i+2$. If $i$ is even, then \eqref{,,6} and \eqref{,,8} imply that
\begin{equation} \label{,,13}
  h(\gamma(t_j)) \pin\le\pin h^+_L(\gamma(t_j)) \pin\le\pin h(\gamma(t_i)) + L\rho(\gamma(t_i),\gamma(t_j)).
\end{equation}
If $i$ is odd these inequalities are true all the more, because in this case $h(\gamma(t_i)) = h^+_L(\gamma(t_i))$, and the function $h^+_L(x)$ is $L$-Lipschitz.

Similarly, using \eqref{,,6} and \eqref{,,10} we obtain the inequalities
\begin{equation} \label{,,14}
  h(\gamma(t_j)) \pin\ge\pin h^-_L(\gamma(t_j)) \pin\ge\pin h(\gamma(t_i)) - L\rho(\gamma(t_i),\gamma(t_j)),
\end{equation}
which, along with \eqref{,,13}, give the Lipschitz condition for the restriction of $h$ to $\gamma(T)$.

Finally, for any $t_i\in T$ and $x'\in X'$ by virtue of \eqref{,,6}, \eqref{,,7} and Corollary \ref{..3} we have
\begin{equation*}
  \bigl| h(\gamma(t_i)) -h(x')\bigr| \le L\rho(\gamma(t_i),x').
\end{equation*}
This proves that the function $h$ is $L$-Lipschitz on the whole $X'\cup \gamma(T)$. \qed

\begin{lemma} \label{..5}
Suppose\/ $h$ is an\/ $L'$-Lipschitz function defined on a finite subset\/ $X'\subset X$, and\/ $\gamma\colon [a,b]\to X$ is an injective continuous mapping. Then for any\/ $L>L'$ and any\/ $\theta\in (0,1)$ there exist a partition\/ $T=\{t_0,t_1,\dots,t_n\}$ of the segment\/ $[a,b]$ and an\/ $L$-Lipschitz extension of\/ $h$ to\/ $X'\cup \gamma(T)$ such that
\begin{equation}\label{,,15}
  \var(h\circ\gamma,T) \pin=\pin \sum_{i=1}^n \big| h(\gamma(t_{i-1})) -h(\gamma(t_i))\big| \pin\ge\pin \theta L\rho(\gamma(a),\gamma(b)).
\end{equation}
\end{lemma}

Below (in Lemma \ref{..8}) the requirement for $\gamma$ to be injective will be removed.

\medskip

\proof. Evidently, the set $\gamma^{-1}(X')$ is finite. Adding to it the points $a,\,b$, we obtain some partition $S = \{s_0,s_1,\dots,s_m\}$ of the segment $[a,b]$. We may choose so long segments $[a_i,b_i]$ in the intervals $(s_{i-1},s_i)$ as to ensure the inequalities
\begin{equation} \label{,,16}
  \sum_{i=1}^m \rho(\gamma(a_i),\gamma(b_i)) \pin\ge\pin \sum_{i=1}^m \sqrt{\theta}\pin \rho(\gamma(s_{i-1}),\gamma(s_i)) \pin\ge\pin \sqrt{\theta}\pin \rho(\gamma(a),\gamma(b)).
\end{equation}
By construction, all the images $\gamma([a_i,b_i])$ do not intersect with each other and with $X'$.

Let us fix numbers
\begin{equation} \label{,,17}
  \max\bigl\{L',\sqrt{\theta}\pin L\bigr\} < L_1 < L_2 < \,\dotsc\, < L_m < L.
\end{equation}
Sequentially applying Lemma \ref{..4}, we construct some partitions $T_i$ of the segments $[a_i,b_i]$ and $L_i$-Lipschitz extensions of $h$ to the sets $X'\cup \gamma(T_1)\cup \,\dotsc\, \cup \gamma(T_i)$ such that
\begin{equation} \label{,,18}
  \var(h\circ\gamma, T_i) \pin\ge\pin L_i\pin \rho(\gamma(a_i), \gamma(b_i)), \qquad i=1,\dots,m.
\end{equation}
Now set $T =\{a\}\cup T_1 \cup \,\dotsc\, \cup T_m \cup \{b\}$. Then summing up \eqref{,,18} and taking into account \eqref{,,16}, \eqref{,,17}, we obtain \eqref{,,15}. \qed

\medskip

Lemma \ref{..5} enables us to prove the sufficient part of Theorem \ref{..1} in the case of injective mapping $\gamma\colon [a,b] \to X$. The proof will be by contradiction. Namely, assume that $\gamma\notin \mathit{AC}([a,b],X)$. Then there exists an $\eps>0$ such that for any natural number $n$ there is a finite set of disjoint intervals $I_n =\{(a_i,b_i)\}$ lying in the segment $[a,b]$ and satisfying the conditions
\begin{equation} \label{,,19}
 \sum_{(a_i,b_i)\in I_n} (b_i -a_i) < \frac{1}{n}, \qquad \sum_{(a_i,b_i)\in I_n} \rho(\gamma(a_i),\gamma(b_i)) >\eps.
\end{equation}

\medskip\noindent
Let us put $I =\bigcup_n I_n$ and introduce a total numbering for all intervals from $I$ by one natural index $i$, so that $I =\{(a_i,b_i)\}_{i\in\mathbb N}$.

Fix an infinite sequence of real numbers
\begin{equation} \label{,,20}
  1 =L_1<L_2<L_3< \dotsc < L = 2.
\end{equation}
Set $T_1=\{a_1,b_1\}$, \,$X'_1 =\gamma(T_1)$ and define an $L_1$-Lipschitz function $h$ on $X'_1$ by the rule
\begin{equation*}
  h(\gamma (a_1)) =0, \qquad h(\gamma(b_1)) =\rho(\gamma(a_1),\gamma(b_1)).
\end{equation*}
Sequentially applying Lemma \ref{..5}, we may construct finite partitions $T_i$ of the segments $[a_i,b_i]$ and $L_i$-Lipschitz extensions of $h$ to the sets $X'_i = \gamma(T_1) \cup \,\dotsc\, \cup \gamma(T_i)$ such that
\begin{equation} \label{,,21}
  \var(h\circ\gamma, T_i) \pin\ge\pin \theta L_i\pin \rho(\gamma(a_i), \gamma(b_i)),
  \qquad i=2,3,\dots
\end{equation}
As a result, we obtain an $L$-Lipschitz function $h$ defined on the set $X' =\bigcup_i X'_i$, which by means of Lemma \ref{..2} extends to an $L$-Lipschitz function on $X$.

Summing up \eqref{,,21} over all $(a_i,b_i)\in I_n$ and taking into account \eqref{,,20}, \eqref{,,19}, we obtain
\begin{equation*}
 \sum_{(a_i,b_i)\in I_n} \var(h\circ\gamma, T_i) \pin\ge\pin
 \sum_{(a_i,b_i)\in I_n} \theta L_i\pin \rho(\gamma(a_i), \gamma(b_i)) \pin>\pin \theta\eps.
\end{equation*}
On the other hand, due to the left inequality in \eqref{,,19}, the total length of all segments of the partitions $T_i$ involved in the last formula is less than $1/n$. Therefore the composition $h\circ\gamma$ is not absolutely continuous, and Theorem \ref{..1} is proved.

\medskip

Note that if Lemma \ref{..5} did not require $\gamma$ to be injective, then the same reasoning would serve as a proof of Theorem \ref{..1} in the general case. To get rid of this injectivity requirement, we will make a suitable piecewise injective modification of $\gamma$.

\smallskip

A \emph{piecewise injective curve} is any continuous mapping $\gamma\colon A\to X$ that has the following properties:

a) it is defined on a compact subset $A$ of the real axis;

b) if $\gamma(c) =\gamma(d)$, where $c,d\in A$, then the interval $(c,d)$ does not intersect with $A$;

c) for each segment $[c,d]$ the set $\gamma([c,d]\cap A)$ is connected.

\smallskip

From b) it follows that the mapping $\gamma$ can take the same values at no more than two different points, and only at the ends of a `hole' in its domain $A$.

\smallskip

A \emph{piecewise injective modification} of the continuous mapping $\gamma\colon [a,b] \to X$ is any piecewise injective curve that is a restriction of $\gamma$ to some compact subset $A\subset [a,b]$ containing the points $a,b$.

\begin{lemma} \label{..6}
Every continuous mapping\/ $\gamma\colon [a,b] \to X$ has a piecewise injective modification.
\end{lemma}

\proof. Let $A_1 =[a,b]$. Then for each $i=1,2,\dots$ we take the longest segment $[c_i,d_i]\subset A_i$ for which $\gamma(c_i) =\gamma(d_i)$ (if there are many of them, then take the leftmost one) and set $A_{i+1} =A_i\setminus (c_i,d_i)$.

Set $A =\bigcap_i A_i$. Then the restriction of $\gamma$ to $A$ is a piecewise injective curve.

Let us check this. First, the set $A$ is compact.

Secondly, assume that $\gamma(c) =\gamma(d)$ for a pair of different points $c,d\in A$. Then take the largest $i\in\mathbb N$ such that $d_i-c_i\ge d-c$. If $[c,d]\ne [c_j,d_j]$ for every $j=1,\dots,i$ then by construction $[c,d]\cap (c_j,d_j) =\varnothing$ and, therefore, $[c,d]\subset A_{i+1}$, which contradicts the choice of~$i$. It follows that the segment $[c,d]$ coincides with one of $[c_j,d_j]$, and then, by construction, the interval $(c,d)$ does not intersect with $A$.

Thirdly, by means of standard topological arguments it is verified that for any segment $[c,d]$ the sets $\gamma([c,d]\cap A_i)$ are connected. It follows that the set $\gamma([c,d]\cap A)$, being the intersection of a sequence of nested connected compact sets $\gamma([c,d]\cap A_i)$, is compact and connected. \qed

\medskip

Now we state and prove an analogue to Lemma \ref{..4} for a piecewise injective curve.

\smallskip

A \emph{partition of a set} $A\subset\mathbb R$ we call any finite collection of points $T =\{t_0,t_1,\dots,t_n\} \subset A$ going in ascending order: $t_0<t_1<\dotsc < t_n$.

\begin{lemma} \label{..7}
Suppose\/ $h$ is an\/ $L'$-Lipschitz function defined on a certain subset\/ $X'\subset X$, and\/ $\gamma\colon A\to X$ is a piecewise injective curve satisfying the condition\/ $\rho(\gamma(A),X') >0$. Then for each\/ $L>L'$ there exist a partition\/ $T=\{t_0,t_1,\dots,t_n\}$ of\/ $A$ and an\/ $L$-Lipschitz extension of the function\/ $h$ to\/ $X'\cup \gamma(T)$ such that
\begin{equation}\label{,,22}
  \var(h\circ\gamma,T) \pin=\pin
  \sum_{i=1}^n \big| h(\gamma(t_{i-1})) -h(\gamma(t_i))\big| \pin\ge\pin L \diam (\gamma(A)).
\end{equation}
\end{lemma}

\proof{\pin} is in essence the same as for Lemma~\ref{..4}. Let $h^{\pm}_L(x)$ be the functions from Lemma~\ref{..2}. Then, by point c) of this lemma, $h^{+}_L(x) > h^{-}_L(x)$ for all $x\in \gamma(A)$.

Take two points $a,b\in A$ such that
\begin{equation}\label{,,23}
  \rho(\gamma(a),\gamma(b)) =\diam (\gamma(A)), \qquad a<b.
\end{equation}

Starting from $t_0 =a$, we will sequentially choose $t_i$ and define alternating values
\begin{equation}\label{,,24}
  h(\gamma(t_{2i})) =h^-_L(\gamma(t_{2i})), \qquad
  h(\gamma(t_{2i+1})) =h^+_L(\gamma(t_{2i+1})),
\end{equation}
except the last point $t_n=b$, where it will be
\begin{equation}\label{,,25}
  h^-_L(\gamma(t_{n})) \pin\le\pin h(\gamma(t_{n})) \pin\le\pin  h^+_L(\gamma(t_{n})),
\end{equation}
according to the next rule: for even $i$ we take
\begin{gather} \label{,,26}
  t_{i+1} \pin=\pin \sup\bigl\{ t\in [t_i,b]\cap A : h^-_L(\gamma(t_i)) +L\rho(\gamma(t_i),\gamma(t)) \le h^+_L(\gamma(t))\bigr\}, \\[6pt] \label{,,27}
  h(\gamma(t_{i+1})) \pin=\pin h^-_L(\gamma(t_i)) +L\rho(\gamma(t_i),\gamma(t_{i+1})),
\end{gather}
and for odd $i$ we take
\begin{gather}  \label{,,28}
  t_{i+1} \pin=\pin \sup\bigl\{ t\in [t_i,b]\cap A : h^+_L(\gamma(t_i)) - L\rho(\gamma(t_i),\gamma(t)) \ge h^-_L(\gamma(t))\bigr\}, \\[6pt]  \label{,,29}
  h(\gamma(t_{i+1})) \pin=\pin h^+_L(\gamma(t_i)) -L\rho(\gamma(t_i),\gamma(t_{i+1})).
\end{gather}

By the definition of a piecewise injective curve, the set $\gamma([t_{i+1},b]\cap A)$ is connected. It follows that in the case $t_{i+1}<b$, depending on the parity of $i$, the inequalities in the right-hand sides of \eqref{,,26}, \eqref{,,28} turns into equalities at the point $t=t_{i+1}$:
\begin{equation*}
 h^\mp_L(\gamma(t_i)) \pm L\rho(\gamma(t_i),\gamma(t_{i+1})) \pin=\pin h^\pm_L(\gamma(t_{i+1})).
\end{equation*}
By virtue of \eqref{,,27} and \eqref{,,29}, they are equivalent to \eqref{,,24}.

\pagebreak[1]

Since the mapping $\gamma\colon A\to X$ is uniformly continuous, and the functions $h^{+}_L(x)$ and $h^{-}_L(x)$ are $L$-Lipschitz, this procedure in a finite number of steps will reach the point $t_n=b$, at which inequalities \eqref{,,25} are satisfied.

From \eqref{,,24}, \eqref{,,27}, \eqref{,,29} we obtain the equalities
\begin{equation} \label{,,30}
  \bigl|h(\gamma(t_{i-1})) -h(\gamma(t_{i}))\bigr| \pin=\pin L\rho(\gamma(t_{i-1}),\gamma(t_i)), \qquad i=1,\dots,n.
\end{equation}
Summing them up and taking into account \eqref{,,23}, we get
\begin{equation*}
  \var(h\circ\gamma,T) \pin=\pin \sum_{i=1}^n L\rho(\gamma(t_{i-1},\gamma(t_i)) \pin\ge\pin L\rho(\gamma(a),\gamma(b)) \pin=\pin L\diam(\gamma(A)).
\end{equation*}

Let us check Lipschitz condition for the restriction of $h$ to $\gamma(T)$. For any pair of adjacent points $\gamma(t_{i-1})$, $\gamma(t_{i})$ it is proved in \eqref{,,30}. Consider an arbitrary pair of points $\gamma(t_{i})$, $\gamma(t_{j})$, where $j\ge i+2$. If $i$ is even, then \eqref{,,24} and \eqref{,,26} imply that
\begin{equation} \label{,,31}
  h(\gamma(t_j)) \pin\le\pin h^+_L(\gamma(t_j)) \pin\le\pin h(\gamma(t_i)) + L\rho(\gamma(t_i),\gamma(t_j)).
\end{equation}
If $i$ is odd these inequalities are true all the more, because in this case $h(\gamma(t_i)) = h^+_L(\gamma(t_i))$, and $h^+_L(x)$ is $L$-Lipschitz. Similarly, using \eqref{,,24} and \eqref{,,28} we obtain the inequalities
\begin{equation} \label{,,32}
  h(\gamma(t_j)) \pin\ge\pin h^-_L(\gamma(t_j)) \pin\ge\pin h(\gamma(t_i)) - L\rho(\gamma(t_i),\gamma(t_j)),
\end{equation}
which, along with \eqref{,,31}, give the Lipschitz condition for the restriction of $h$ to $\gamma(T)$.

Finally, for any $t_i\in T$ and $x'\in X'$ by virtue of \eqref{,,24}, \eqref{,,25} and Corollary \ref{..3} we have
\begin{equation*}
  \bigl| h(\gamma(t_i)) -h(x')\bigr| \le L\rho(\gamma(t_i),x').
\end{equation*}
This proves that the function $h$ is $L$-Lipschitz on the whole $X'\cup \gamma(T)$. \qed

\begin{lemma} \label{..8}
Lemma \ref{..5} is valid for any continuous mapping\/ $\gamma\colon [a,b]\to X$ $($i.\,e., without the injectivity condition\/$)$.
\end{lemma}

\proof. Take a piecewise injective modification $\gamma\colon A\to X$ such that $a,b\in A$. Then the intersection $A\cap\gamma^{-1}(X')$ is finite. Adding to it the points $a,\,b$, we obtain some partition $S = \{s_0,s_1,\dots,s_m\}$ of $A$. Exploiting the connectivity of the sets $\gamma([s_{i-1},s_i]\cap A)$, we may choose points $a_i,b_i\in (s_{i-1},s_i)\cap A$ such that
\begin{equation} \label{,,33}
  \sum_{i=1}^m \rho(\gamma(a_i),\gamma(b_i)) \pin\ge\pin \sum_{i=1}^m \sqrt{\theta}\pin \rho(\gamma(s_{i-1}),\gamma(s_i)) \pin\ge\pin \sqrt{\theta}\pin \rho(\gamma(a),\gamma(b)).
\end{equation}
By construction, all images $\gamma([a_i,b_i]\cap A)$ do not intersect with each other and with $X'$.

Let us fix numbers
\begin{equation} \label{,,34}
  \max\bigl\{L',\sqrt{\theta}\pin L\bigr\} < L_1 < L_2 < \,\dotsc\, < L_m < L.
\end{equation}
Applying Lemma \ref{..7} to the piecewise injective curves $\gamma\colon [a_i,b_i]\cap A \to X$, we sequentially construct some partitions $T_i$ of the sets $[a_i,b_i]\cap A$ and some $L_i$-Lipschitz extensions of $h$ to the sets $X'\cup \gamma(T_1)\cup \,\dotsc\, \cup \gamma(T_i)$ such that
\begin{equation} \label{,,35}
  \var(h\circ\gamma, T_i) \pin\ge\pin L_i\pin \rho(\gamma(a_i), \gamma(b_i)), \qquad i=1,\dots,m.
\end{equation}
Define the partition $T =\{a\}\cup T_1 \cup \,\dotsc\, \cup T_m \cup \{b\}$ of the segment $[a,b]$. Summing up \eqref{,,35} and taking into account \eqref{,,33}, \eqref{,,34}, we obtain \eqref{,,15}. \qed

\medskip

The sufficient part of Theorem \ref{..1} in the general case is proved in exactly the same way as in the case of the injective mapping $\gamma$, provided Lemma \ref{..8} is used instead of Lemma \ref{..5}.

\end{document}